\documentclass[11pt]{article} 
\usepackage{amssymb}

   \csname @addtoreset\endcsname{equation}{section}
    \csname @addtoreset\endcsname{figure}{section}
    \csname @addtoreset\endcsname{table}{section}

\newcounter{thanksnum}
\def\thanksnumber#1
{\setcounter{thanksnum}{\value{footnote}}\setcounter{footnote}{#1}%
                     \addtocounter{footnote}{-1}\footnotemark
                     \setcounter{footnote}{\value{thanksnum}}}
\def\newtheoremz#1{\@ifnextchar[{\@othmz{#1}}{\@nthmz{#1}}}

\def\@nthmz#1#2{%
\@ifnextchar[{\@xnthmz{#1}{#2}}{\@ynthmz{#1}{#2}}}

\def\@xnthmz#1#2[#3]{\expandafter\@ifdefinable\csname #1\endcsname
{\@definecounter{#1}\@addtoreset{#1}{#3}%
\expandafter\xdef\csname the#1\endcsname{\expandafter\noexpand
  \csname the#3\endcsname \@thmcountersepz \@thmcounterz{#1}}%
\global\@namedef{#1}{\@thmz{#1}{#2}}\global\@namedef{end#1}{\@endtheoremz}}}

\def\@ynthmz#1#2{\expandafter\@ifdefinable\csname #1\endcsname
{\@definecounter{#1}%
\expandafter\xdef\csname the#1\endcsname{\@thmcounterz{#1}}%
\global\@namedef{#1}{\@thm{#1}{#2}}\global\@namedef{end#1}{\@endtheoremz}}}

\def\@othmz#1[#2]#3{\expandafter\@ifdefinable\csname #1\endcsname
  {\global\@namedef{the#1}{\@nameuse{the#2}}%
\global\@namedef{#1}{\@thmz{#2}{#3}}%
\global\@namedef{end#1}{\@endtheoremz}}}

\def\@thmz#1#2{\refstepcounter
    {#1}\@ifnextchar[{\@ythmz{#1}{#2}}{\@xthmz{#1}{#2}}}

\def\@xthmz#1#2{\@begintheoremz{#2}{\csname the#1\endcsname}\ignorespaces}
\def\@ythmz#1#2[#3]{\@opargbegintheoremz{#2}{\csname
       the#1\endcsname}{#3}\ignorespaces}

\def\@thmcounterz#1{\noexpand\arabic{#1}}
\def\@thmcountersepz{.}
\def\@begintheoremz#1#2{ \trivlist \item[\hskip \labelsep{\bf #1\ #2}]}
\def\@opargbegintheoremz#1#2#3{ \trivlist
      \item[\hskip \labelsep{\bf #1\ #2\ (#3)}]}
\def\@endtheoremz{\endtrivlist}

\newtheorem{theorem}{Theorem}[section]
\newtheorem{proposition}[theorem]{Proposition}
\newtheorem{lemma}[theorem]{Lemma}

\newtheorem{definition}[theorem]{Definition}
\newtheorem{remark}[theorem]{Remark}

\newtheorem{condition}[theorem]{Condition}

\def\defi{\stackrel{{\scriptscriptstyle \Delta}}{=}}

\def\d{\delta}
\def\o{\omega}
\def\O{\Omega}

\def\F{{\cal F}}
\def\w{\widehat}

\def\esssup{\mathop{\rm ess\, sup}}
\def\const{{\rm const\,}}

\def\R{{\bf R}}
\def\E{{\bf E}}
\def\P{{\bf P}}

\def\L{L}

\def\b{\beta}
\def\s{\delta}
\def\g{\gamma}

\def\W{{\cal W}^*}
\def\ww{\widetilde}

\def\t{\theta}
\def\oo{\bar}
\def\s{\sigma}

\def\p{\partial}

\def\U{{\cal U}}

\def\A{{\cal A}}
\def\M{{\cal M}}

\def\L{{\cal L}}
\def\I{{\, \cal I}}

\newcommand{\be}{\begin{equation}}
\newcommand{\ee}{\end{equation}}
\newcommand{\bd}{\begin{displaymath}}
\newcommand{\ed}{\end{displaymath}}
\newcommand{\ba}{\begin{array}{ll}}
\newcommand{\ea}{\end{array}}
\newcommand{\baa}{\begin{eqnarray}}
\newcommand{\eaa}{\end{eqnarray}}
\newcommand{\baaa}{\begin{eqnarray*}}
\newcommand{\eaaa}{\end{eqnarray*}}   \font\sm=cmr10

\def\ww{\tilde}
\def\W{{\cal W}}

\def\Q{{\cal Q}}
\def\CC{{\cal C}}

\def\LU{L}

\def\KK{{K}}

\def\I{{ \d}}

\title{Duality and semi-group property for backward parabolic Ito
equations \footnote{{\em Random Operators and Stochastic Equations.
} (2010) {\bf 18}, 51-72.}}
\author{ Nikolai Dokuchaev
\\ {\sm  Department of Mathematics and Statistics,
Curtin University of Technology},\\ {\sm GPO Box U1987, Perth,
Western Australia, 6845}}
\begin{document}
\maketitle
\begin{abstract}
We study  existence, uniqueness, semi-group property, and a priori
estimates for solutions for backward parabolic Ito equations in
domains with boundary. We study also  duality between forward and
backward equations. The semi-group for backward equations is
established in the form of some anti-causality.  The novelty is that
 the semi-group property involves the diffusion term that is a part
of the solution.
{\it Key
words and phrases:} parabolic Ito equations, backward SPDEs,
regularity, semi-group property, anti-causality
\end{abstract}

\section{Introduction}
The paper studies backward stochastic partial differential equations
(SPDEs) with a Dirichlet boundary condition in a cylinder
$D\times[0,T]$ for a region $D\subseteq \R^n$ with boundary
condition at the terminal time $t=T$.   The difference between
backward and forward equations is not that important for the
deterministic equations since a backward equation can be converted
to a forward equation by a time change. However, it cannot be done
so easily for stochastic equations, because one looks for  a
solution adapted to the driving Brownian motion. It is why the
backward stochastic differential equations
 require special
consideration.  A possible approach is to consider the so-called
Bismut backward equations such the diffusion term is not given a
priori but needs to be found. These backward SPDEs were widely
studied (see, e.g., Yong and Zhou
 (1999) and references there;
 non-linear ordinary backward Ito equations was studied in Pardoux and Peng
 (1990);
backward SPDEs was studied also by Dokuchaev (1992),
 (2003)). Note that there is a duality between linear forward and
backward equations.   Forward SPDEs were studied in the literature
(see, e.g., Al\'os et al (1999), Bally {\it et al} (1994),
Chojnowska-Michalik and Goldys (1995), Da Prato and Tubaro (1996),
Gy\"ongy (1998), Krylov (1999), Maslowski (1995), Pardoux (1993),
 Rozovskii (1990), Walsh (1986), Zhou (1992), Dokuchaev
 (1995), (2002), (2005),
and the bibliography there). Duality between forward and backward
equations was studied by Zhou (1992) for domains without boundary
and by Dokuchaev (1992) for domains with boundary in some special
cases. Backward SPDEs represent analogs of backward parabolic
Kolmogorov equations for non-Markov Ito processes, including the
case of bounded domains, so they may be used for characterization of
distributions of first exit time in non-Markovian setting, as was
shown by the author (1992).  A different type of backward equations
was described in Chapter 5 of Rozovskii (1990).

\par In the present paper, we study
existence, uniqueness, a priori estimates, duality, and semi-group
properties for solutions for backward linear parabolic Ito equations
in domains. The novelty is that we consider domains with boundaries;
the semi-group property involves the diffusion term that is a part
of the solution. The proofs for prior estimates are based on duality
between forward and backward equations, so the most part of the
paper is devoted to establishing this duality. This duality  is used
also to establish that backward parabolic equations have some
causality (more precisely, some anti-causality). This fact can help,
for instance, to split time intervals and apply numerical methods
via using of auxiliary problems with pathwise constant in time
coefficients, or to apply dynamic programming methods for the
corresponding control problems.
\section{Definitions}\label{SecD}
\subsection{Spaces and classes of functions.} 
Assume that we a given an open domain $D\subseteq\R^n$ such that
either $D=\R^n$ or $D$ is bounded  with $C^2$-smooth boundary $\p
D$. Let $T>0$ be given, and let $Q\defi D\times (0,T)$. \par We
are given a standard complete probability space $(\O,\F,\P)$ and a
right-continuous filtration $\F_t$ of complete $\s$-algebras of
events, $t\ge 0$. We are given also a $N$-dimensional process
$w(t)=(w_1(t),...,w_N(t))$ with independent components such that
it is a Wiener process with respect to $\F_t$.
\par
We denote by $\|\cdot\|_{ X}$ the norm in a linear normed space
$X$, and
 $(\cdot, \cdot )_{ X}$ denote  the scalar product in  a Hilbert space $
X$.
\par
We introduce some spaces of real valued functions.
\par
 Let $G\subset \R^k$ be an open
domain, then ${W_q^m}(G)$ denote  the Sobolev  space of functions
that belong to $L_q(G)$ with the distributional derivatives up to
the $m$th order,  $q\ge 1$.
\par
 We denote Euclidean norm in $\R^k$ as $|\cdot|$, and $\bar G$ denote
the closure of a region $G\subset\R^k$.
\par Let $H^0\defi L_2(D)$,
and let $H^1\defi \stackrel{\scriptscriptstyle 0}{W_2^1}(D)$ be
the closure in the ${W}_2^1(D)$-norm of the set of all smooth
functions $u:D\to\R$ such that  $u|_{\p D}\equiv 0$. Let
$H^2=W^2_2(D)\cap H^1$ be the space equipped with the norm of
$W_2^2(D)$. The spaces $H^k$ and $W_2^k(D)$ are called  Sobolev
spaces, they are Hilbert spaces, and $H^k$ is a closed subspace of
$W_2^k(D)$, $k=0,1,2$.
\par
 Let $H^{-1}$ be the dual space to $H^{1}$, with the
norm $\| \,\cdot\,\| _{H^{-1}}$ such that if $u \in H^{0}$ then
$\| u\|_{ H^{-1}}$ is the supremum of $(u,v)_{H^0}$ over all $v
\in H^1$ such that $\| v\|_{H^1} \le 1 $. $H^{-1}$ is a Hilbert
space.
\par We shall write $(u,v)_{H^0}$ for $u\in H^{-1}$
and $v\in H^1$, meaning the obvious extension of the bilinear form
from $u\in H^{0}$ and $v\in H^1$.
\par
We denote by $\oo\ell _{k}$ the Lebesgue measure in $\R^k$, and we
denote by $ \oo{{\cal B}}_{k}$ the $\sigma$-algebra of Lebesgue
sets in $\R^k$.
\par
We denote by $\oo{{\cal P}}$  the completion (with respect to the
measure $\oo\ell_1\times\P$) of the $\s$-algebra of subsets of
$[0,T]\times\O$, generated by functions that are progressively
measurable with respect to $\F_t$.
\par
Let $Q_s\defi D\times [s,T]$. For $k=-1,0,1,2$, we  introduce the
spaces
 \baaa
 X^{k}(s,T)\defi L^{2}\bigl([ s,T ]\times\Omega,
{\oo{\cal P }},\oo\ell_{1}\times\P;  H^{k}\bigr),\quad Z^k_t
\defi L^2\bigl(\Omega,{\cal F}_t,\P;
H^k\bigr),\quad \CC^{k}(s,T)\defi C\left([s,T]; Z^k_T\right),
\eaaa
Furthermore, introduce the spaces $$ Y^{k}(s,T)\defi
X^{k}(s,T)\!\cap \CC^{k-1}(s,T), \quad k\ge 0, $$ with the norm $
\| u\| _{Y^k(s,T)}
\defi \| u\| _{{X}^k(s,T)} +\| u\| _{\CC^{k-1}(s,T)}. $
\par
For brevity, we shall use the notations
 $X^k\defi X^k(0,T)$, $\CC^k\defi \CC^k(0,T)$,
and  $Y^k\defi Y^k(0,T)$.
\par
  The
spaces $X^k$ and $Z_t^k$  are Hilbert spaces.
\par
 Further, introduce the spaces
\baaa \W^{k}_p \defi L^{2}\bigl([0,T ]\times\O, \overline{{\cal
P}},\oo\ell_{1}\times\P;\, W_p^{k}(D)\bigr), \quad k=0,1,\ldots,
\quad 1\le p\le
 +\infty.
\eaaa The same notation we shall use the space of $n\times
n$-dimensional matrix functions. In that case,
$\|\cdot\|_{\W^k_p}$ means the summa of all this norms for all
components.
\begin{proposition} 
\label{propL} Let $\xi\in X^0$,
 let a sequence  $\{\xi_k\}_{k=1}^{+\infty}\subset
L^{\infty}([0,T]\times\O, \ell_1\times\P;\,C(D))$ be such that all
$\xi_k(\cdot,t,\o)$ are progressively measurable with respect to
$\F_t$, and let $\|\xi-\xi_k\|_{X^0}\to 0$. Let $t\in [0,T]$ and
$j\in\{1,\ldots, N\}$ be given.
 Then the sequence of
the integrals $\int_0^t\xi_k(x,s,\o)\,dw_j(s)$ converges in
$Z_t^0$ as $k\to\infty$, and its limit depends on $\xi$, but does
not depend on $\{\xi_k\}$.
\end{proposition}
\par
{\it Proof} follows from completeness of  $X^0$ and from the
equality
\begin{eqnarray*}
\E\int_0^t\|\xi_{k}(\cdot,s,\o)-\xi_m(\cdot,s,\o)\|_{H^0}^2\,ds
=\int_D\,dx\,\E\left(\int_0^t\big(\xi_k(x,s,\o)-
\xi_m(x,s,\o)\big)\,dw_j(s)\right)^2.
\end{eqnarray*}
\begin{definition} 
\rm Let $\xi\in X^0$, $t\in [0,T]$, $j\in\{1,\ldots, N\}$, then we
define $\int_0^t\xi(x,s,\o)\,dw_j(s)$ as the limit  in $Z_t^0$ as
$k\to\infty$ of a sequence $\int_0^t\xi_k(x,s,\o)\,dw_j(s)$, where
the sequence $\{\xi_k\}$ is such  as in Proposition \ref{propL}.
\end{definition}
 Sometimes we shall omit
$\o$.
\par
For $t\in[0,T]$, define operators $\I_t: C([0,T];Z_T^{k})\to
Z^{k}_t$ such that $\I_tu=u(\cdot,t)$, where $k=-1,0,1$.
\section{Review of existence results for forward equations}
\label{SecC} Let $s\in [0,T)$, $\varphi\in X^{-1}$, $h_i\in X^0$,
and $\Phi\in Z^0_s$. Consider the problem \be \label{4.1} \ba
d_tu=\left( \A u+ \varphi\right)dt + \sum_{i=1}^N[
B_iu+h_i]dw_i(t),\quad \quad t\ge s,\\ u|_{t=s}=\Phi,\quad
u(x,t,\o)|_{x\in \p D}=0 . \ea
 \ee
 Here
 $u=u(x,t,\o)$,
 $(x,t)\in Q$,   $\o\in\O$, and
  \be\label{A}\A v\defi \sum_{i,j=1}^n
b_{ij}(x,t,\o)\frac{\p^2v}{\p x_i \p x_j}(x) +\sum_{i=1}^n
f_{i}(x,t,\o)\frac{\p v}{\p x_i }(x) +\,\lambda(x,t,\o)v(x), \ee
where $b_{ij}, f_i, x_i$ are the components of $b$,$f$, and $x$.
Further, \be\label{B} B_iv\defi\frac{dv}{dx}\,(x)\,\beta_i(x,t,\o)
+\oo \beta_i(x,t,\o)\,v(x),\quad i=1,\ldots ,N. \ee
\par
We assume that the functions $b(x,t,\o):
\R^n\times[0,T]\times\O\to\R^{n\times n}$, $\b_j(x,t,\o):
\R^n\times[0,T]\times\O\to\R^n$, $\oo\b_i(x,t,\o):$
$\R^n\times[0,T]\times\O\to\R$, $f(x,t,\o):
\R^n\times[0,T]\times\O\to\R^n$, $\lambda(x,t,\o):
\R^n\times[0,T]\times\O\to\R$ and  $\varphi (x,t,\o): \R^n\times
[0,T]\times\O\to\R$ are  progressively measurable for any $x\in
\R^n$ with respect to $\F_t$.
\par We assume that $b(x,t,\o)$, $f(x,t,\o)$, $\lambda(x,t,\o)$ vanish for
$(x,t,\o)\notin D\times [0,T]\times\Omega $.\par
 To proceed further, we assume that Conditions
\ref{cond3.1.A}-\ref{condB} remain in force throughout this paper.
 \begin{condition} \label{cond3.1.A}  (Coercivity) The matrix  $b=b^\top$ is
symmetric,  bounded, and progressively measurable with respect to
$\F_t$ for all $x$, and there exists a constant $\d>0$ such that
\be
 \label{Main1} y^\top  b
(x,t,\o)\,y-\frac{1}{2}\sum_{i=1}^N |y^\top\b_i(x,t,\o)|^2 \ge
\d|y|^2 \quad\forall\, y\in \R^n,\ (x,t)\in  D\times [0,T],\
\o\in\O. \ee
\end{condition}
\par
Inequality (\ref{Main1})  means that equation (\ref{4.1}) is
coercive or {\it superparabolic}, in the terminology of Rozovskii
(1990).
\begin{condition}\label{cond3.1.B}
Functions $b(x,t,\o):\R^n \times \R\times\O\to \R^{n \times n}$,
$f(x,t,\o):\R^n \times \R\times\O\to \R^n$, $\lambda (x,t,\o):\R^n
\times \R\times\O\to \R$,  are bounded and differentiable in $x$,
and \baaa \esssup_{\o}\esssup_{(x,t)\in    Q} \biggl[ \Bigl|
\frac{\p b}{\p x}(x,t,\o)\Bigr| + \Bigl| \frac{\p f}{\p
x}(x,t,\o)\Bigr| + \Bigl| \frac{\p \lambda}{\p x}(x,t,\o)\Bigr|
\biggl]< +\infty. \eaaa
\end{condition}
\begin{condition}\label{condB}
The functions  $\b_i(x,t,\o)$ and $\oo\b_i(x,t,\o)$ are bounded
and differentiable in $x$, and  $\esssup_{x,t,\o}|\frac{\p
\b_i}{\p x}(x,t,\o)|<+\infty$, $\esssup_{x,t,\o}|\frac{\p\oo
\b_i}{\p x}(x,t,\o)|<+\infty$, $i=1,\ldots ,N$.
\end{condition}
\par
We introduce the set of parameters   $$ \ba {\cal P}
\defi \biggl( n,\,\, D,\,\,  T\,\, \delta,\,\,\,\,
\esssup_{x,t,\o}\Bigl[| b(x,t,\o)|+ |f(x,t,\o)|+ \Bigl|\frac{\p
b}{\p x}(x,t,\o)\Bigr|+ \Bigl|\frac{\p f}{\p
x}(x,t,\o)\Bigr|\Bigr],\\
 \esssup_{x,t,\o,i}\Bigl[|
\b_i(x,t,\o)|+ |\oo\b_i(x,t,\o)|+ \Bigl|\frac{\p \b_i}{\p
x}(x,t,\o)\Bigr|+ \Bigl|\frac{\p\oo\b_i}{\p x}(x,t,\o)\Bigr|\Bigr]
\biggr). \ea $$
\subsubsection*{The definition of solution} 
\begin{definition} 
\label{defsolltion} \rm Let $h_i\in X^0$ and $\varphi\in X^{-1}$.
We say that   equations (\ref{4.1}) are satisfied for $u\in Y^1$
if
\begin{eqnarray}
&&u(\cdot,t,\o)-u(\cdot,r,\o)\nonumber
\\  &&\hphantom{xxx}= \int_r^t\big(\A u(\cdot,s,\o)+
\varphi(\cdot,s,\o)\big)\,ds+ \sum_{i=1}^N
\int_r^t[B_iu(\cdot,s,\o)+h_i(\cdot,s,\o)]\,dw_i(s)
\label{intur}
\end{eqnarray}
for all $r,t$ such that $0\le r<t\le T$, and this equality is
satisfied as an equality in $Z_T^{-1}$.
\end{definition}
Note that the condition on $\p D$ is satisfied in the following
sense:   $u(\cdot,t,\o)\in H^1$ for a.e. \ $t,\o$.  Further, $u\in
Y^1$, and  the value of  $u(\cdot,t,\o)$ is uniquely defined in
$Z_T^0$ given $t$, by the definitions of the corresponding spaces.
The integrals with $dw_i$ in (\ref{intur}) are defined as elements
of $Z_T^0$. The integral with $ds$ is defined as an element of
$Z_T^{-1}$. Definition \ref{defsolltion} requires for (\ref{4.1})
that this integral must be equal  to an element of $Z_T^{0}$ in
the sense of equality in $Z_T^{-1}$.
\subsubsection*{Existence theorems and fundamental
inequalities for forward equations}
\par The following Lemma combines the first and the second fundamental inequalities and related existence result
for forward SPDEs (the cases when $k=-1$ and $k=1$ respectively). It
is an analog of the so-called "energy inequalities", or "the
fundamental inequalities" known for deterministic parabolic
equations (Ladyzhenskaya {\it et al} (1969)).
\begin{lemma}
\label{lemma1} Let either $k=-1$ or $k=0$. Assume that Conditions
\ref{cond3.1.A}, \ref{cond3.1.B},  and \ref{condB}, are satisfied.
In addition,  assume that if $k=0$, then $\b_i(x,t,\o)=0$ for
$x\in \p D$, $i=1,...,N$. Let $\varphi\in X^{k}(s,T)$ and $\Phi\in
Z_s^{k+1}$. Then problem (\ref{4.1}) has an unique solution $u$ in
the class $Y^1(s,T)$, and the following analog of the first
fundamental inequality is satisfied:  \be \label{4.2} \| u
\|_{Y^{k+2}(s,T)}\le c \left(\| \varphi \|
_{X^{k}(s,T)}+\|\Phi\|_{Z^{k+1}_s}+ \sum_{i=1}^N\|h_i
\|_{X^{k+1}(s,T)}\right), \ee where  $c=c({\cal P})$ is a constant
that depends on ${\cal P} $ only.
\end{lemma}
\par
The result  of Lemma  \ref{lemma1}  for $k=-1$ is well known for
long time (see, e.g., Rozovskii (1990), Ch. 3.4.1).  The result
for $k=0$ was obtained  in Dokuchaev (2005).
\par
Note that $Y^{k+2}(s,T)=X^{k+2}(s,T)\!\cap \CC^{k+1}(s,T)$, hence
the solution $u=u(\cdot,t)$ is continuous in $t$ in
$L_2(\O,\F,\P,H^{k+1})$.
\par
 Introduce  operators $L(s,T):X^{-1}(s,T)\to
Y^1(s,T)$, $\M_{i}(s,T):X^{0}(s,T)\to Y^1(s,T)$, and
$\L(s,T):Z^0_s\to Y^1(s,T)$, such that
$$u=L(s,T)\varphi+\L(s,T)\Phi+\sum_{i=1}^N\M_{i}(s,T)h_i,$$ where
$u$ is the solution in $Y^1(s,T)$ of   problem (\ref{4.1}). These
operators are linear and continuous; it follows immediately from
Lemma \ref{lemma1}. We shall denote by $L$, $\M_i$, and $\L$, the
operators $L(0,T)$, $\M_i(0,T)$, and $\L(0,T)$, correspondingly.
\section{Backward  equations}\label{secBE} Introduce the operators being formally adjoint to
the operators $\A$ and $B_i$:
\begin{eqnarray*}
\A^*v &=&\sum_{i,j=1}^n\frac{\p ^2 }{\p x_i \p x_j}\,
\biggl(b_{ij}(x,t,\o)\,v(x)\biggr)
-\, \sum_{i=1}^n \frac{\p }{\p x_i }\,\big(
f_{i}(x,t,\o)\,v(x)\big) +\lambda(x,t,\o)\,v(x),
\\
 B_i^*v&=&-\sum_{i=1}^n \frac{\p }{\p x_i }\,\big(\beta_{i}(x,t,\o)\,v(x))+\oo \beta_i(x,t,\o)\,v(x).
\end{eqnarray*}
Consider the boundary value problem in $Q$
\begin{eqnarray}
\begin{array}{c}
\label{4.8}  d_tp+
\left(\A^*p+\sum_{i=1}^NB_i^*\chi_i+\xi\right)\,dt=
\sum_{i=1}^N\chi_i \,dw_i(t),
\\
\noalign{\vspace*{4pt}} \quad p|_{t=T}=\Psi, \quad
p(x,t,\o)\,|_{x\in \p D}=0.
\end{array}
\end{eqnarray}
\begin{definition} 
\label{defsolltion2} \rm We say that equation (\ref{4.8}) is
satisfied for $p\in Y^1$, $\Psi\in Z_T^0$, $\chi_i\in X^0$ if
\begin{eqnarray} 
\label{intur1} p(\cdot,t)=\Psi+\int_t^T\Biggl(\A^*p(\cdot,s) +
\sum_{i=1}^NB_i^*\chi_i(\cdot,s)+\xi(\cdot,s)\Biggr)\,ds
-\sum_{i=1}^N \int_t^T\chi_i(\cdot,s)\,dw_i(s)
\end{eqnarray}
for any $t\in[0,T]$. The equality here is assumed to be an
equality in the space $Z_T^{-1}$.
\end{definition}
\begin{theorem} 
\label{Th4.2} For any $\xi \in X^{-1}$ and $\Psi\in Z_T^0$, there
exists a pair $(p,\chi)$, such that $p\in Y^1$,
$\chi=(\chi_1,\ldots, \chi_N)$, $\chi_i\in X^0$ and (\ref{4.8}) is
satisfied. This pair is uniquely defined, and the following analog
of the first fundamental inequality is satisfied: \be\label{1FE}
\|p\|_{Y^1}+\sum_{i=1}^N\|\chi_i\|_{X^0}\le
c(\|\xi\|_{X^{-1}}+\|\Psi\|_{Z_T^{0}}), \ee where $ c=c({\cal
P})>0$ is a constant that does not depend on $\xi$ and $\Psi$. In
addition, $$p={\LU}^*\xi+(\I_T\LU)^*\Psi,\quad
\chi_i=\M_i^*\xi+(\I_T\M_i)^*\Psi,\quad
p(\cdot,0)=\L^*\xi+(\I_T\L)^*\Psi,$$ where ${\LU}^*: X^{-1}\to
X^1$, ${\M}_i^*: X^{0}\to X^0$, $(\I_T\LU)^*: Z_0^{0}\to X^1$,
$(\I_T\M_i)^*: Z_0^{0}\to X^0$, and $(\I_T\L)^*: Z_T^{0}\to
Z_0^0$,
 are the operators that are adjoint  to the operators ${\LU}:
X^{-1}\to X^1$, $\M_i : X^{0}\to X^1$, $\I_T\M_i: X^{-1}\to Z_T^0$,
$\I_T\M_i: X^0\to Z_T^{0}$, and $\,\I_T\L: Z_0^0\to Z_T^{0}$,
respectively.
\end{theorem}
An example of application of duality established in this theorem is
given in Theorem \ref{lemmaMark}.
\section{ Proof of Theorem \ref{Th4.2}}
Theorem \ref{Th4.2}
 extends Theorem 4.1 from Dokuchaev (1992) and
 Theorem 3.2  in Dokuchaev (2003) for the case of
 non-zero $\Psi,\b_i,\oo \b_i$. Apparently, the extension
 required was non-trivial and, as can be seen from the proof below, required solid efforts.
We shall prove Theorem \ref{Th4.2} using the following steps:
first, we obtain some  decomposition results for the basic
operators, the the proof for the case when the coefficients of
$\A$ are $\F_0$-measurable and $B_i\equiv 0$, then we consider the
case when $\A$ is of the general form and $B_i\equiv 0$, and then
we consider the general case.
\subsection{Decomposition of operators $L$ and
$\M_i$}\label{SecDecomposition} Our method of  proof is based on
decomposition of the operators to superpositions of simpler
operators.
\begin{definition} 
\label{Def4.2} \rm Define operators $\KK: Z_0^{0}\to Y^1$, $\Q_0:
X^{-1}\to Y^1$, $\Q_i:\allowbreak X^0\to Y^1$, $i=1,...,N$, as the
operators $\L: Z_0^{0}\to Y^1$, $\LU: X^{-1}\to Y^1$,
$\M_i:\allowbreak X^0\to Y^1$, $i=1,...,N$, considered for the case
when $B_i=0$ for all $i$.
\end{definition}
By Lemma \ref{lemma1}, these linear operators are continuous. It
follows from the definitions  that
$$
\KK\Phi+\Q_0\eta+\sum_{i=1}^N\Q_ih_i=V,
$$
where $\eta\in X^{-1}$, $\Phi\in Z_0^0$, and $h_i\in X^0$, and
where $V$ is the solution of the problem
\begin{eqnarray}
\begin{array}{c}
\label{4.11} d_tV=(\A V+\eta)\,dt + \sum_{i=1}^N h_i\,dw_i(t),
\\
 V|_{t=0}=\Phi,\quad V(x,t,\o)\,|_{x\in \p D}=0.
\end{array}
\end{eqnarray}
 Define the operators \be\label{P*} P\defi\sum_{i=1}^N\Q_iB_i,
\quad P^*\defi\sum_{i=1}^NB_i^*\Q_i^*. \ee Clearly, the operator
$P: X^{1}\to X^{1}$ is continuous, and $P^*: X^{-1}\to X^{-1}$ is
its adjoint  operator. Hence the operator $P^*: X^{-1}\to X^{-1}$
is continuous. Since the operators $\Q_i^*:
\left(\CC^{0}\right)^*\to X^0$ and $B_i^*: X^0\to X^{-1}$ are
continuous, it follows that the operator  $P^*:
\left(\CC^{0}\right)^*\to X^{-1}$ is continuous.
\par
Let $$ P_0\defi \I_T\sum_{i=1}^N\Q_iB_i,\quad P_0^*\defi
\sum_{i=1}^NB_i^*(\I_T\Q_i)^*. $$
\begin{lemma} 
\label{prop4.3} The operator $(I-P)^{-1}: X^{1}\to X^{1}$ is
continuous, and \baa &&{\LU}=(I-P)^{-1}\Q_0, \quad\L=(I-P)^{-1}K,
\quad \M_i=(I-P)^{-1}\Q_i,\nonumber \\
&&\I_T\LU=P_0(I-P)^{-1}\Q_0+\I_T\Q_0,\quad
\I_T\M_i=P_0(I-P)^{-1}\Q_i+\I_T\Q_i, \label{I-P}\eaa $i=1,...,N$.
The operator $(I-P^*)^{-1}: X^{-1}\to X^{-1}$ is also continuous,
and \baa &&\LU^*=\Q_0^*(I-P^*)^{-1},\quad
\L^*=\KK^*(I-P^*)^{-1},\quad \M_i^*=\Q_i^*(I-P^*)^{-1},\nonumber
\\
&&(\I_T\LU)^*=\Q_0^*(I-P^*)^{-1}P_0^*+(\I_T\Q_0)^*, \quad
(\I_T\M_i)^*=\Q_i^*(I-P^*)^{-1}P_0^*+(\I_T\Q_i)^*.\hphantom{xxxx}
\label{I-P*} \eaa \end{lemma}
\par {\it Proof}. If
$\g\in X^{1}$, then the solution of  problem (\ref{4.1}) with
$\Phi=0$, $\varphi=0$, and $h_i\defi B_i\g$,
 can be represented as $$
 V=\sum_{i=1}^n\Q_iB_i(\g +V)=P(\g+V).
 $$
 By Lemma  \ref{lemma1},
$\|V\|_{Y^1}\le\const\|\g\|_{X^1}$. Hence the problem
\begin{eqnarray}
\label{zP} V=P(\g+V)
\end{eqnarray}
has the unique solution $V\in Y^{1}$. Define the operator $\w
P:X^{1}\to Y^{1}$ such that $V=\w P\g= P(\g +V)$ is the solution of
problem (\ref{zP}).
\par
Let $\g\in X^{1}$ be given. Set  $V\defi \w P\g$, $h\defi\g+V$.
Clearly, $\|h\|_{X^1}\le\const(\|\g\|_{X^1}+\|\zeta\|_{X^1})?$. By
the definitions, $$V=Ph,\quad \g=h-Ph,\quad h=(I-P)^{-1}\g. $$ We
have noticed already that  $\|V\|_{Y^1}\le\const\|\g\|_{X^1}$. Hence
$\|h\|_{X^1}\le\const\|\g\|_{X^1}$, and the operator $(I-P)^{-1}:
X^{1}\to X^{1}$ is continuous. Therefore, the adjoint operator
$(I-P^*)^{-1}: X^{-1}\to X^{-1}$ is also continuous.
\par
By the definitions,  the solution of problem (\ref{4.2}) has the
form $$ V={\LU}\,\varphi+\L\Phi+\sum_{i=1}^N\M_ih_i=\Q_0\varphi +
\sum_{i=1}^n\Q_i[B_iV+h_i]+\KK\Phi =\Q_0\varphi
+\sum_{i=1}^N\Q_ih_i+PV+\KK\Phi$$ for any $\varphi\in X^{-1}$,
$h_i\in X^0$, and $\Phi\in Z_0^0$. Then the first three equations
in (\ref{I-P}) follow. Clearly,  \baaa (I-P)^{-1}=P(I-P)^{-1}-I.
\eaaa
 Hence the last two equations in (\ref{I-P}) follow.  Therefore,
all equations (\ref{I-P*}) for the adjoint operators  follow. This
completes the proof of Lemma \ref{prop4.3}. $\Box$

\par
We found  that the operator $(I-P^*)^{-1}:X^{-1}\to X^{-1}$ is
continuous. In addition, Theorem 3.2 from p. 467 of Dokuchaev
(1992) stated that the operator $P^*: X^0\to X^0$ is continuous
under some additional technical conditions. The question arises if
the operator $(I-P)^{-1}:X^0\to X^0$ is continuous. If this is
true, then, by Lemma \ref{prop4.3}, the operators $L^*:X^0\to
Y^2$, $\M_i^*:X^0\to X^1$ are continuous; in that case, an analog
of the second fundamental inequality holds for backward equations,
and existence has place for solutions  $p\in X^2$ and $\chi_i\in
X^1$.

\subsection{Proof of Theorem \ref{Th4.2} for $B_i\equiv 0$ and
 $\F_0$-measurable $(b,f,\lambda)$ }
\par
 The following lemma  extends Theorem 4.1 from Dokuchaev (1992)  for the case of non-zero $\Psi$.
\begin{lemma}\label{lemmaLapl} Theorem \ref{Th4.2} holds under additional assumptions
that  the function $\mu=(b,f,\lambda)$ is such that $\mu(x,t,\o)$ is
$\F_0$-measurable for all $(x,t)$, and $\b_i\equiv 0$,
$\oo\b_i\equiv 0$ for $i=1,...,N$.
\end{lemma}
\par
{\it Proof}. It suffices to prove that,  for any $\xi \in X^{-1}$,
there exists a pair $(p,\chi)$, with $p\in Y^1$,  and
$\chi=(\chi_1,\ldots, \chi_N)$, $\chi_i\in X^0$, such that
\begin{eqnarray}
\begin{array}{c}
\label{4.14} d_tp+(\A^*p+\xi)\,dt=\sum_{i=1}^N\chi_i\,dw_i(t),
\\
\noalign{\vspace*{5pt}} \quad p|_{t=T}=\Psi,\quad p(x,t,\o)\,|_{x\in
\p D}=0,
\end{array}
\end{eqnarray}
and that \be p=\Q_0^*\xi+(\I_T\Q_0)^*\Psi,\quad
\chi_i=\Q_i^*\xi+(\I_T\Q_i)^*\Psi, \quad
p(.,0)=\KK^*\xi+(\I_T\KK)^*\Psi, \label{pQ}\ee where $\Q_0^*:
X^{-1}\to X^1$, $\Q_i^*: X^{-1}\to X^0$, $(\I_T\Q_0)^*:  Z^0_T\to
X^1$, $(\I_T\Q_i)*: Z_T^{0}\to X^0$, $i=1,...,N$, and $\KK^*:
X^{-1}\to Z_0^0$, are the operators that are adjoint  for the
operators $\Q_0: X^{-1}\to X^1$, $\Q_i: X^{0}\to X^1$, $\I_T\Q_i:
X^{0}\to Z_T^0$, and  $\KK^*:  Z_0^0\to X^1$, correspondingly. We
know that these operators are continuous.
\par We interpret the solution of problem (\ref{4.14}) similarly
to Definition \ref{defsolltion2}.
\par
Let $\oo p$ be the solution of the problem \be\ba \frac{\p \oo p}{\p
t}+\A^*\oo p=-\xi,\quad\\ \oo p|_{t=T}=\Psi,\quad \oo
p(x,t,\o)|_{x\in \p D}=0. \label{oop}\ea\ee  Clearly, $\oo p\in \oo
X^1\cap C([0,T];Z_T^0)$.
\par
Let $p\defi{\cal E}\oo p$, and let $\varphi\in
L_{\infty}(\O,\P,\F_T,C(Q))$, $\Phi\in L_2(\O,\P,\F_0,C^1(D))\cap
Z_0^1$. Then $\oo p$ and its derivatives presented in (\ref{oop})
belong to $L_{\infty}(\O,\P,\F_T,C(Q))$. We have that if
$V=\Q_0\varphi+\KK\Phi$ then \be\ba \frac{\p V}{\p t}=\A
V+\varphi,\quad\\ V|_{t=0}=\Phi,\quad V(x,t,\o)|_{x\in \p D}=0.
\label{ooV}\ea\ee Hence $$ (\oo p(\cdot,T),V(\cdot,T))_{Z_T^0}-(\oo
p(\cdot,0),V(\cdot,0))_{Z_T^0} =(-\A^*\oo p-\xi,V)_{\oo X^0}+(\oo
p,\A V+\varphi)_{\oo X^0}. $$
$$ (\Psi,V(\cdot,T))_{Z_T^0}-(\oo p(\cdot,0),\Phi)_{Z_T^0}
=(-\xi,V)_{\oo X^0}+(\oo p,\varphi)_{\oo X^0}. $$ Hence \be
(\Psi,V(\cdot,T))_{Z_T^0}-(p(\cdot,0),\Phi)_{Z_0^0} =(-\xi,V)_{
X^0}+(p,\varphi)_{X^0}, \label{dual1}\ee i.e., \baaa
(p,\varphi)_{X^0}+(p(\cdot,0),\Phi)_{Z_0^0}&=&(\Psi,V(\cdot,T))_{Z_T^0}+(\xi,V)_{X^0}\\
&=&(\Psi,\I_T \Q_0\varphi+\I_T\KK
\Phi)_{Z_T^0}+(\xi,\Q_0\varphi+\KK\Phi)_{X^0}. \eaaa Then
$p=\Q_0^*\xi+(\I_T\Q_0)^*\Psi$ and
$p(.,0)=\KK^*\xi+(\I_T\KK)^*\Psi$.
\begin{remark}\label{remAA}
Up to this point, we didn't use the assumption that the coefficients
$\A$ are $\F_0$-measurable; all previous reasons are valid for the
general $\A$ as well.
\end{remark}
\par
Further, by Clark's Theorem, there exist functions
$\g_i(\cdot,t,\cdot)\in X^0$, $\g_{\Psi i}\in X^0$, and $\g_{\xi
i}(\cdot,t,\cdot)\in X^0$, such that \baaa \oo p(x,t,\o)=\E\oo
p(x,t,\o)+\sum_{i=1}^N\int_0^T\g_i(x,t,s,\o)dw_i(s),\\\oo
\Psi(x,t,\o)=\E\Psi(x,\o)+\sum_{i=1}^N\int_0^T\g_{\Psi
i}(x,s,\o)dw_i(s),\\
\xi(x,t,\o)=\E\xi(x,t,\o)+\sum_{i=1}^N\int_0^T\g_{\xi
i}(x,t,s,\o)dw_i(s). \eaaa Moreover, it follows that ${\cal D}
g_i(\cdot,t,\cdot)\in X^0$, where either ${\cal D}\g=\p\g/\p t$ or
${\cal D} \g=\A^*\g=\Delta\g$, and \baaa {\cal D}\oo
p(x,t,\o)=\E{\cal D}\oo p(x,t,\o)+\sum_{i=1}^N\int_0^T{\cal
D}\g_i(x,t,s,\o)dw_i(s). \eaaa By (\ref{oop}), \be\ba \frac{\p
\g_i}{\p t}(\cdot,t,s,\o)+\A^*\g_i(\cdot,t,s,\o)=-\g_{\xi
i}(\cdot,t,s,\o),\quad\\ \g_i(x,T,s,\o) =\g_{\Psi i}(x,s,\o),\quad
\g_i(x,t,s,\o)|_{x\in \p D}=0. \label{oog}\ea\ee It follows that
\baaa \sup_{t\in[s,T]}\int_D|\g_i(x,t,s,\o)|^2dx\le
c\left[\int_s^T\|\g_{\xi
i}(\cdot,t,s,\o)\|_{H^{-1}}^2dt+\int_D|\g_{\Psi
i}(x,s,\o)|^2dx\right], \eaaa where $c=c(T,n,D)>0$ is a constant.
Hence \baaa \int_D|\g_i(x,s,s,\o)|^2dx\le c\left[\int_s^T\|\g_{\xi
i}(\cdot,t,s,\o)\|_{H^{-1}}^2dt+\int_D|\g_{\Psi
i}(x,s,\o)|^2dx\right]. \eaaa Let \be \label{chi1}
\chi_i(x,t,\o)\defi \g_i(x,t,t,\o).\ee It follows that $$
\|\chi_i\|_{X^0}\le c(\|\xi\|_{X^{-1}}+ \|\Psi\|_{Z_T^0}). $$
\par
Let us show that the pair $(p,\chi)$, $\chi=(\chi_1,\ldots,
\chi_N)$, is such that (\ref{4.14}) is satisfied. Clearly, \baaa \oo
p(x,t,\o)= p(x,t,\o)+\sum_{i=1}^N\int_t^T\g_i(x,t,s,\o)dw_i(s),
\eaaa and \baaa \oo p(\cdot,t)=\Psi+\int_t^T\Bigl(\A^*\oo p(\cdot,s)
+\xi(\cdot,s)\Bigr)\,ds.\eaaa Hence \baaa
&&p(\cdot,t)=\Psi+\int_t^T\Bigl(\A^*p(\cdot,s)
+\xi(\cdot,s)\Bigr)\,ds \\
&&+\sum_{i=1}^N\biggl[-\int_t^T
\g_i(\cdot,t,s)dw_i(s)+\int_t^Tds\int_s^T[
\A^*\g_i(\cdot,s,r)+\g_{\xi i}(\cdot,s,r)]dw_i(r)\biggr]
\\
&&=\Psi+\int_t^T\Bigl(\A^*p(\cdot,s) +\xi(\cdot,s)\Bigr)\,ds \\
&&+\sum_{i=1}^N\biggl[-\int_t^T \g_i(\cdot,t,s)dw_i(s)+\int_t^T
dw_i(r)\int_t^r[
\A^*\g_i(\cdot,s,r)+\g_{\xi i}(\cdot,s,r)]ds\biggr]\\
\\
&&=\Psi+\int_t^T\Bigl(\A^*p(\cdot,s) +\xi(\cdot,s)\Bigr)\,ds
\\&&+\sum_{i=1}^N\int_t^T dw_i(s)\biggl[-\g_i(\cdot,t,s)+\int_t^s[
\A^*\g_i(\cdot,r,s)+\g_{\xi
i}(\cdot,r,s)]dr\biggr]\\
&&=\Psi+\int_t^T\Biggl(\A^*p(\cdot,s)+\xi(\cdot,s)\Biggr)\,ds
-\,\int_t^T\sum_{i=1}^N \chi_i(\cdot,s)\,dw_i(s), \eaaa
\par
since
$$
\g_i(\cdot,t,s)-\int_t^s[ \A^*\g_i(\cdot,r,s)+\g_{\xi
i}(\cdot,r,s)]dr =\g_i(\cdot,s,s)=\chi_i(\cdot,s).
$$
 We have that if
$V=\Q_ih$, where $h\in X^0$, then \be\ba d_tV=\A Vdt
+hdw_i(t),\quad\\ V|_{t=0}=0,\quad V(x,t,\o)|_{x\in \p D}=0.
\label{Vh}\ea\ee Hence $$ (
p(\cdot,T),V(\cdot,T))_{Z_T^0}-(p(\cdot,0),V(\cdot,0))_{Z_T^0}
=(-\A^*\oo p-\xi,V)_{ X^0}+( p,\A V)_{X^0}+(\chi_i,h)_{X^0}, $$
i.e., $ (\Psi,V(\cdot,T))_{Z_T^0}=(-\xi,V)_{X^0}+(\chi_i,h)_{X^0}. $
It follows that $$
(\chi_i,h)_{X^0}=(\Psi,\I_T\Q_ih)_{Z_T^0}+(\xi,\Q_ih)_{X^0} \quad
\forall h\in X^0. $$ Hence function (\ref{chi1}) is \be \label{chi2}
\chi_i=\Q_i^*\xi+(\I_T\Q_i)^*\Psi. \ee This completes the proof of
Lemma \ref{lemmaLapl}.
\subsection{Proof of Theorem \ref{4.2} for $B_i\equiv 0$}
\begin{lemma}\label{lemmaB=0} Theorem \ref{4.2} holds under additional assumptions
that  $\b_i\equiv 0$, $\oo\b_i\equiv 0$ for $i=1,...,N$, i.e., when
$B_i= 0$ for all $i$.
\end{lemma}
\par
{\it Proof.} It suffices to prove that (\ref{pQ}) holds. Introduce
the operators $\ww\A$ and $\ww\A^*$ such that $$ \ww \A v=\ww
\A^*\defi \Delta v=\sum_{i=1}^n\frac{\p ^2v}{\p x_i^2}, $$ where
$\Delta$ is the Laplace operator,Denote by $\ww \Q_i$, $\ww\KK$,
$\ww \Q^*_i$ and $\ww\KK^*$,  the operators, defined similarly
$\Q_i$,
 $\KK$,  $\Q_i^*$, $\KK^*$,  but with substituting $\A=\ww \A$. Introduce the operators
$$ \U\defi\ww\Q_0(\ww \A- \A), \quad \U_{\Delta}\defi\ww\Q_0^*(\ww
\A^*-\A^*), \quad .
$$
Clearly, the operators $\A: X^{1}\to X^{-1}$ and $\A^*: X^{1}\to
X^{-1}$ are continuous. By Lemma \ref{lemma1}, the operator
 $\Q_0: X^{-1}\to X^1$ is continuous. Hence the adjoint
 operator  $\Q_0^*: X^{-1}\to X^1$ is continuous, the operators ${\U:
X^{1}\to X^{1}}$, $\U_{\Delta}: X^{1}\to X^{1}$ are both continuous,
and the adjoint  operators $\U^*: X^{-1}\to X^{-1}$ and
$\U_{\Delta}^*: X^{-1}\to X^{-1}$ are continuous, where
$\U_{\Delta}^*=(\ww A- \A)\,\ww\Q_0$ and $\U^*\defi(\ww A^*-
\A^*)\,\ww\Q_0^*$.

1. Let us prove that the operator $(I+\U_{\Delta})^{-1}\ww\Q_0^*:
X^{1}\to X^{1}$ is continuous. First, we shall prove that the
operator $(I+\U)^{-1}\ww\Q_0: {X^{-1}\to X^{-1}}$ is continuous.

Let $\eta\in X^{-1}$, $h_i\in X^{0}$, and $\Phi\in Z_0^0$. By Lemma
\ref{lemma1}, the boundary value problem \baaa \label{z} &&d_tz
=\big(\ww Az+( A-\ww A)\,z+\eta\big)\,dt +\sum_{i=1}^N h_i\,dw_i(t)
=(Az+\eta)\,dt+\sum_{i=1}^N h_i\,dw_i(t),
\nonumber\\
&&z|_{t=0}=\Phi, \quad z(x,t,\o)\,|_{x\in \p D}=0 \label{3.7} \eaaa
has the unique solution $z\in Y^1$, and
$\|z\|_{Y^1}\le\const(\|\eta\|_{X^{-1}}+\sum_i\|h_i\|_{X^0}+\|\Phi\|_{Z_0^{0}})$.
By the definitions of the corresponding operators, the solution $z$
of the equation (\ref{z}) is
$$z=\Q_0\eta+\sum_i\Q_ih_i+\KK\Phi =\ww\Q_0(\eta +(\A-\ww \A)\,z)+
\sum_i\ww\Q_ih_i+\ww\KK\Phi,$$ or
\begin{eqnarray}
\label{z1} z=\Q_0\eta+\sum_{i=1}^N\Q_ih_i+\KK\Phi
=(I+\U)^{-1}\Bigg(\ww\Q_0\eta
+\sum_{i=1}^N\ww\Q_ih_i+\ww\KK\Phi\Bigg). \label{3.8}
\end{eqnarray}
Hence the operator $(I+\U)^{-1}\ww\Q_0: X^{-1}\to X^{1}$ is
continuous.
\par
Further, note that the range of the operator $\ww\Q_0(X^{-1})$
contains all smooth functions from $X^{1}$, therefore, the range is
dense in $X^{-1}$. Clearly, the equality $\ww\Q_0y=0$ implies that
$y=0$. Let $y\in \ww\Q_0(X^{-1})$ and $y=\ww\Q_0z$. Then the
equality $\ww\Q_0(I+\U_{\Delta}^*)^{-1}\,x=y$ implies equalities
\baaa (I+\U_{\Delta}^*)^{-1}x=z, \quad (I+\U_{\Delta}^*)\,z=x,\quad
\ww\Q_0^{-1}y+\U_{\Delta}^*\ww\Q_0^{-1}y=x,\quad
\ww\Q_0^{-1}y+(\ww A-A)\,y=x,\\
\ww\Q_0^{-1}y+\ww\Q_0^{-1}\U\,y=x,\quad
\ww\Q_0^{-1}(I+\U)\,y=x,\eaaa i.e.\
\begin{eqnarray}
\label{q1} \ww\Q_0(I+\U_{\Delta}^*)^{-1}=(I+\U)^{-1}\,\ww\Q_0.
\end{eqnarray}
This and continuity of the operator $(I+\U)^{-1}\ww\Q_0: X^{-1}\to
X^{1}$ imply that the operator $\ww\Q_0(I+\U_{\Delta}^*)^{-1}:
X^{-1}\to X^{1}$ is continuous. Therefore, the adjoint  operator
$(I+\U_{\Delta})^{-1}\ww\Q_0^*: X^{-1}\to X^{1}$ is continuous.
\par
2. Let us establish the connection between (\ref{4.14}) and the
operators $\Q_i$, ${i=0,\ldots, N}$. It was shown above that the
operator $(I+\U_{\Delta})^{-1}\ww\Q_0^*: X^{-1}\to X^{1}$ is
continuous. By (\ref{z1}), it follows that $$
\Q_0\xi=(I+\U)^{-1}\,\ww\Q_0\xi,\quad
\Q_ih_i=(I+\U)^{-1}\ww\Q_ih_i,\quad i=1,...,N,\quad
\KK\Phi=(I+\U)^{-1}\,\ww\KK\Phi
$$ for all $\eta\in X^{-1}$,
$h_i\in X^{0}$, and $\Phi\in Z_0^0$. Hence
 \be\label{q0qi}
\Q_0=(I+\U)^{-1}\,\ww\Q_0,\quad \Q_i=(I+\U)^{-1}\ww\Q_i,\quad
i=1,...,N,\quad \KK=(I+\U)^{-1}\,\ww\KK, \ee  and \be\label{q0q2}
\Q_0^*=\ww\Q_0^*(I+\U^*)^{-1},\quad
\Q_i^*=\ww\Q_i^*(I+\U^*)^{-1},\quad i=1,...,N,\quad
\KK=\,\ww\KK^*(I+\U^*)^{-1}. \ee By (\ref{q1}), it follows that
\begin{eqnarray}
\label{q2} (I+\U_{\Delta})^{-1}\ww\Q_0^*=\ww\Q_0^*(I+\U^*)^{-1}.
\end{eqnarray}
\par
Set $$ p\defi(I+\U_{\Delta})^{-1}[\ww\Q_0^*\xi+(\I_T\ww
\Q_0)^*\Psi]. $$ Clearly, $(I+\U_{\Delta})\,p=\ww\Q_0^*\xi+(\I_T\ww
\Q_0)^*\Psi$, i.e. $$ p=\ww\Q_0^*\xi-\U_{\Delta}p +(\I_T\ww
\Q_0)^*\Psi=\ww\Q_0^*(\xi+(\A^*-\ww \A^*)\,p)+(\I_T\ww \Q_0)^*\Psi.
$$  For this $p$, set $$ \chi_i\defi\ww\Q_i^*(\xi+(\A^*-\ww
\A^*)\,p)+(\I_T\ww \Q_i)^*\Psi.$$ Let us use Lemma \ref{lemmaLapl}
now. By this lemma and and by the definitions for the corresponding
operators, $(p,\chi_1,\ldots, \chi_N)\in Y^1\times (X^0)^N$
satisfies (\ref{4.14}). By (\ref{q2}) and (\ref{q0qi}), it follows
that $$ p=\ww\Q_0^*(I+\U^*)^{-1}\xi +(I+\U_{\Delta})^{-1}(\I_T\ww
\Q_0)^*\Psi =\Q_0^*\xi +(\I_T \Q_0)^*\Psi. $$ We have used here that
$(\I_T \Q_0)^*=\Q_0^*\I_T^*$, where $\Q_0^*:
\left(\CC^{0}\right)^*\to X^1$ is the linear continuous operator
that is adjoint  to the linear continuous operator $\Q_0: X^{-1}\to
\CC^0$, and $\I_T^*: Z_T^0\to \left(\CC^{0}\right)^*$ is the linear
continuous operator that is adjoint  to the linear continuous
operator $\I_T^*: \CC^0\to Z_T^0$. Here $\left(\CC^{0}\right)^*$ is
the space that is dual for $\CC^0$.
\par
 Further, $\chi_i$ satisfies
\baaa \chi_i &=&\ww\Q_i^*\big(\xi+(\A^*-\ww \A^*)\,[\Q_0^*\xi+(\I_T
\Q_0)^*\Psi]\big)+(\I_T\ww \Q_i)^*\Psi
\\ &=&\ww\Q_i^*\big(\xi+(\A^*-\ww
\A^*)\,\ww\Q_0^*(I+\U^*)^{-1}\xi\big)+ \ww\Q_i^*(\A^*-\ww
\A^*)\,(\I_T \Q_0)^*\Psi+(\I_T\ww \Q_i)^*\Psi
\\
&=&\ww\Q_i^*\big(I-\U^*(I+\U^*)^{-1}\big)\xi+\ww\Q_i^*\big(I-\U^*(I+\U^*)^{-1}\big)\I_T^*\Psi\\&=&
\ww\Q_i^*(I+\U^*)^{-1}\,\xi
+\ww\Q_i^*(I+\U^*)^{-1}\I_T^*\Psi\\&=&\Q_i^*\xi+(\I_T\Q_i)^*\Psi.
\eaaa We have used  (\ref{q0qi}) for the last equality. In addition,
we have used again that $(\I_T \Q_i)^*=\Q_i^*\I_T^*$, where $\Q_i^*:
\left(\CC^{0}\right)^*\to X^0$ is the linear continuous operator
that is adjoint  to the linear continuous operator $\Q_0: X^{0}\to
\CC^0$, and $\I_T^*: Z_T^0\to \left(\CC^{0}\right)^*$ is the linear
continuous operator that is adjoint  to the linear continuous
operator $\I_T: \CC^0\to Z_T^0$.
\par
Similarly, \baaa p(\cdot,0) &=&\ww\KK^*\big(\xi+(\A^*-\ww
\A^*)p]\big)+(\I_T\ww \KK)^*\Psi
\\ &=&\ww\KK^*\big(\xi+(\A^*-\ww
\A^*)\,[\Q_0^*\xi+(\I_T \Q_0)^*\Psi]\big)+(\I_T\ww \KK)^*\Psi
\\ &=&\ww\KK^*\big(\xi+(\A^*-\ww
\A^*)\,\ww\Q_0^*(I+\U^*)^{-1}\xi\big)+ \ww\KK^*(\A^*-\ww
\A^*)\,(\I_T \Q_0)^*\Psi+(\I_T\ww \KK)^*\Psi
\\
&=&\ww\KK^*\big(I-\U^*(I+\U^*)^{-1}\big)\xi+\ww\KK_i^*\big(I-\U^*(I+\U^*)^{-1}\big)\I_T^*\Psi\\&=&
\ww\KK^*(I+\U^*)^{-1}\,\xi
+\ww\KK^*(I+\U^*)^{-1}\I_T^*\Psi\\&=&\KK^*\xi+(\I_T\KK)^*\Psi. \eaaa
This completes the proof of Lemma \ref{lemmaB=0}. $\Box$
\begin{remark} Typically, we use the notations such as  $(\I_T L)^*:Z_T^0\to
X^0$ when the entire adjoint operator $\I_TL$ has some sense, so we
don't need to take into account properties of the operator $\I_T^*$.
We use these superposition-type notations for these operators
because it reduces the number of symbols for different operators.
The only exception is the proof of Lemma \ref{lemmaB=0} above, when
we consider operators $(\I_T\Q_i)^*:Z_T^0\to X^0$ as the
superposition of the  operators $\Q_i^*: \left(\CC^{0}\right)^*\to
X^0$ and $\I_T^*: Z_T^0\to \left(\CC^{0}\right)^*$.
\end{remark}
\subsection{Proof of Theorem \ref{Th4.2} for the general case}
Now we are able to prove Theorem \ref{Th4.2}. First, Lemma
\ref{prop4.3} implies that there exists $\zeta\in X^{-1}$ such
that $\zeta=P^*(\xi+\I_T^*\Psi+\zeta)$, i.e.,
$\zeta=(I-P^*)^{-1}P^*(\xi+\I_T^*\Psi)$. Since the operators $P^*:
\left(\CC^{0}\right)^*\to X^{-1}$  and $\I_T^*: Z_T^0\to
\left(\CC^{0}\right)^*$ are continuous, we have that $\zeta\in
X^{-1}$. Set $$p\defi \Q_0^*(\xi+\zeta)+(\I_T\Q_0)^*\Psi,\quad
\chi_i\defi \Q_i^*(\xi+\zeta)+ \Q_i^*\I_T^*\Psi,\quad
\chi\defi(\chi_1,\ldots, \chi_N).$$ Clearly, $(p,\chi)\in
Y^1\times (X^0)^N$. Let us show that the pair $(p,\chi)$ is a
solution of problem (\ref{4.8}). By the definition of the operator
$\Q_i^*$ and by Lemma \ref{lemmaB=0}, it suffices to show that
$\zeta=\sum_{i=1}^NB_i^*\chi_i$.
\par
Let $\oo\zeta\defi\sum_{i=1}^NB_i^*\chi_i$. The equations for
$\chi_i$ are such that
$$
\oo\zeta=\sum_{i=1}^NB_i^*\chi_i=\sum_{i=1}^NB_i^*\Q_i^*(\xi+\zeta)
+\sum_{i=1}^NB_i^* \Q_i^*\I_T^*\Psi=P^*(\xi+\I_T^*\Psi+\zeta).$$ The
definition of  $\zeta$ implies that
 $\oo\zeta=\zeta$. Therefore, the pair $(p,\chi)$ is a solution of  problem (\ref{4.8}).
Farther, we have that \baaa p= \Q_0^*(\xi+\zeta)+(\I_T\Q_0)^*\Psi=
 \Q_0^*(\xi+(I-P^*)^{-1}P^*(\xi+\I_T^*\Psi))+(\I_T\Q_0)^*\Psi
 \\
=
 \Q_0^*(I+(I-P^*)^{-1}P^*)[\xi+(\I_T\Q_0)^*\Psi]
 =\Q_0^*(I-P^*)^{-1}[\xi+(\I_T\Q_0)^*\Psi].
\eaaa (we treat $\I_T$ similarly to the proof of Lemma
\ref{lemmaB=0}). By Lemma \ref{prop4.3}, it follows that
$p=L^*\xi+(\I_T L)^*\Psi$. Similarly, \baaa &&\chi_i=
\Q_i^*(\xi+\zeta)+(\I_T\Q_i)^*\Psi=
 \Q_i^*(\xi+(I-P^*)^{-1}P^*(\xi+\I_T^*\Psi))+(\I_T\Q_i)^*\Psi
 \\
&&=
 \Q_i^*(I+(I-P^*)^{-1}P^*)[\xi+(\I_T\Q_i)^*\Psi]
 =\Q_i^*(I-P^*)^{-1}[\xi+(\I_T\Q_i)^*\Psi]\\&&=\M_i^*\xi+(\I_T
 \M_i)^*\Psi,
\eaaa and \baaa &&p(\cdot,0)= \KK^*(\xi+\zeta)+(\I_T\KK)^*\Psi=
 \KK^*(\xi+(I-P^*)^{-1}P^*(\xi+\I_T^*\Psi))+(\I_T\KK)^*\Psi
 \\
&&=
 \KK^*(I+(I-P^*)^{-1}P^*)[\xi+(\I_T\Q_0)^*\Psi]
 =\KK^*(I-P^*)^{-1}[\xi+(\I_T\Q_0)^*\Psi]\\&&=\KK^*\xi+(\I_T
 \KK)^*\Psi.
 \eaaa
 This completes the proof of Theorem \ref{Th4.2}. $\Box$
\section{Semi-group property for backward equations} It is known that the
dynamic of forward parabolic Ito equation has semi-group property
(causality): if $u=L\varphi+\L_0\Phi$, where $\varphi\in X^{-1}$,
$\Phi\in Z_0^0$, then \be\label{casF}
u|_{t\in[\t,s]}=(L\varphi+\L_0\Phi)|_{t\in[\t,s]}=L(\t,s)\varphi+\L_{\t}(\t,s)u(\cdot,\t).\ee
We establish a similar property for the backward equations
(anti-causality). This property involves the diffusion term that is
a part of the solution.
\begin{theorem}\label{lemmaMark} (Semi-group property for backward equations).
Let $0\le \t< s<T$, and let $p=L^*\xi$, $\chi_i=\M_i\xi$, where
$\xi\in X^{-1}$ and $\Psi\in Z_T^0$. Then \baa
&&p|_{t\in[\t,s]}=L(\t,s)^*\xi|_{t\in[\t,s]}+(\I_sL(\t,s))^*p(\cdot,s),\label{p-cas}
\\
&&p(\cdot,\t)=(\I_s\L_{\t}(\t,s))^*p(\cdot,s)+\L_{\t}(\t,s)^*\xi,\label{p-cas0}
\\
&&\chi_k|_{t\in[\t,s]}=\M_k(\t,s)^*\xi|_{t\in[\t,s]}+(\I_s\M_i(\t,s))^*p(\cdot,s),
\quad k=1,...,N. \label{chi-cas} \eaa
\end{theorem}
\par
{\it Proof of Theorem \ref{lemmaMark}.} Let
$u=L(\t,s)\varphi+\L_{\t}(\t,s)\Phi$, where $\varphi\in X^0(\t,s)$
and $\Phi\in Z_s^1$ are arbitrary. We have that \baaa
&&(p(\cdot,s),u(\cdot,s))_{Z_T^0}-(p(\cdot,\t),\Phi)_{Z_T^0}
\\
&&=\Bigl(-\A^* p-\xi-\sum_{i=1}^NB_i^*\chi_i,u\Bigr)_{
X^0(\t,s)}+(p,\A u+\varphi)_{
X^0(\t,s)}+\sum_{i=1}^N\bigl(\chi_iB_iu\bigr)_{ X^0(\t,s)}. \eaaa
Hence \baa
(p(\cdot,s),u(\cdot,s))_{Z_T^0}-(p(\cdot,\t),\Phi)_{Z_T^0}
=-\bigl(\xi,u\bigr)_{\oo X^0(\t,s)}+(p,\varphi)_{X^0(\t,s)}.
\label{duality}\eaa
 i.e., \baaa
&&(p,\varphi)_{ X^0(\t,s)}+
(p(\cdot,\t),\Phi)_{Z_T^0}=(p(\cdot,s),u(\cdot,s))_{Z_T^0}+\bigl(\xi,u\bigr)_{
X^0(\t,s)}\\ &&=(p(\cdot,s),\I_s[L(\t,s)\varphi+\L_\t(\t,s)\Phi]
)_{Z_T^0} +\bigl(\xi,L(\t,s)\varphi+\L_\t(\t,s)\Phi\bigr)_{
X^0(\t,s)}. \eaaa Hence
 \baaa
 (p,\varphi)_{\oo X^0(\t,s)}+(p(\cdot,\t),\Phi)_{Z_T^0}=((\I_sL(\t,s))^*p(\cdot,s),\varphi)_{ X^0(\t,s)}
+(L(\t,s)^*\xi,\varphi)_{X^0(\t,s)}\\+((\I_s\L_{\t}(\t,s))^*p(\cdot,s),\Phi)_{Z_T^0}
+\bigl(\L_{\t}^*\xi,\Phi\bigr)_{Z_T^0}\\
\hphantom{xxx}=((\I_sL(\t,s))^*p(\cdot,s)+L(\t,s)^*\xi,\varphi)_{
X^0(\t,s)}+((\I_s\L_{\t}(\t,s))^*p(\cdot,s)+\L_{\t}(\t,s)^*\xi,\Phi)_{Z_T^0}.
\eaaa Take $\Phi=0$, then  desired equation (\ref{p-cas})  for
$p|_{[\t,s]}$ follows. Take $\varphi=0$, then  desired equation
(\ref{p-cas0})  for $p(\cdot,\t)$ follows. \par Similarly, let
$u=\M_k(\t,s) h$, where $h\in X^1(\t,s)$ and $k\in\{1,...,N\}$. We
have that \baaa
&&(p(\cdot,s),u(\cdot,s))_{Z_T^0}=(p(\cdot,s),u(\cdot,s))_{Z_T^0}-(p(\cdot,\t),u(\cdot,\t)_{Z_T^0}
\\
&&=\Bigl(-\A^* p-\xi-\sum_{i=1}^NB_i^*\chi_i,u\Bigr)_{
X^0(\t,s)}+(p,\A u)_{
X^0(\t,s)}+\sum_{i=1}^N\bigl(\chi_iB_iu\bigr)_{
X^0(\t,s)}+\bigl(\chi_k,h\bigr)_{ X^0(\t,s)}. \eaaa Hence $
(p(\cdot,s),u(\cdot,s))_{Z_T^0} =-\bigl(\xi,u\bigr)_{\oo
X^0(\t,s)}+(\chi_k,h)_{X^0(\t,s)}, $
 i.e., \baaa (\chi_k,h)_{X^0(\t,s)}&=&(p(\cdot,s),u(\cdot,s))_{Z_T^0}+\bigl(\xi,u\bigr)_{ X^0(\t,s)}\\
&=&(p(\cdot,s),\I_s\M_k(\t,s)h)_{Z_T^0} +\bigl(\xi,\M_k(\t,s)
h\bigr)_{ X^0(\t,s)}. \eaaa Hence $
(\chi_k,h)_{X^0(\t,s)}=((\I_s\M_k(\t,s))^*p(\cdot,s),h)_{Z_T^0}
+\bigl(\M_k(\t,s)^*\xi,h\bigr)_{ X^0(\t,s)} . $ Then desired
equation (\ref{chi-cas}) for $\chi_k$ follows.
 $\Box$

\subsection*{Acknowledgment}  This work  was supported by NSERC
grant of Canada 341796-2008 to the author.
\section*{References}
$\hphantom{XX}$Al\'os, E., Le\'on, J.A., Nualart, D. (1999).
 Stochastic heat equation with random coefficients
 {\it
Probability Theory and Related Fields} {\bf 115}, 1, 41-94.
\par
Bally, V., Gyongy, I., Pardoux, E. (1994). White noise driven
parabolic SPDEs with measurable drift. {\it Journal of Functional
Analysis} {\bf 120}, 484 - 510.
\par
Chojnowska-Michalik, A., and Goldys, B. (1965). {Existence,
uniqueness and invariant measures for stochastic semilinear
equations in Hilbert spaces},  {\it Probability Theory and Related
Fields},  {\bf 102}, No. 3, 331--356.
\par
Da Prato, G., and Tubaro, L. (1996). { Fully nonlinear stochastic
partial differential equations}, {\it SIAM Journal on Mathematical
Analysis} {\bf 27}, No. 1, 40--55.\
\par
Dokuchaev, N.G. (1992). { Boundary value problems for functionals
of
 Ito processes,} {\it Theory of Probability and its Applications}
 {\bf 36 }, 459-476.
\par
 Dokuchaev, N.G. (1995). Probability distributions  of  Ito's
processes: estimations for density functions and for conditional
expectations of integral functionals. {\it Theory of Probability
and its Applications} {\bf 39} (4),  662-670.
\par Dokuchaev, N.G. (2003). Nonlinear
parabolic Ito's equations and duality approach, {\it Theory of
Probability and its Applications} {\bf 48} (1), 45-62.
\par
Dokuchaev, N.G. (2005). Parabolic Ito equations and second
fundamental inequality.  {\it Stochastics} {\bf 77}, iss. 4.,
349-370.
\par
Gy\"ongy, I. (1998). Existence and uniqueness results for
semilinear stochastic partial differential equations. {\it
Stochastic Processes and their Applications} {\bf 73} (2),
271-299.
\par
Kim, Kyeong-Hun (2004). On stochastic partial di!erential
equations with variable coefficients in $C^1$-domains. Stochastic
Processes and their Applications {\bf 112}, 261--283.
\par Krylov, N. V. (1999). An
analytic approach to SPDEs. Stochastic partial differential
equations: six perspectives, 185--242, Math. Surveys Monogr., 64,
Amer. Math. Soc., Providence, RI.
\par
Ladyzhenskaia, O.A. (1985). {\it The Boundary Value Problems of
Mathematical Physics}. New York: Springer-Verlag.
\par
Maslowski, B. (1995). { Stability of semilinear equations with
boundary and pointwise noise}, {\it Ann. Scuola Norm. Sup. Pisa
Cl. Sci.} (4), {\bf  22}, No. 1, 55--93.
\par
Pardoux, E., S. Peng, S. (1990). Adapted solution of a backward
stochastic differential equation. {\it System \& Control Letters}
{\bf 14}, 55-61.
\par
Pardoux, E. (1993).
 Stochastic partial differential equations, a review, {\it Bull. Sc. Math.}
 {\bf 117}, 29-47.
\par
Rozovskii, B.L. (1990). {\it Stochastic Evolution Systems; Linear
Theory and Applications to Non-Linear Filtering.} Kluwer Academic
Publishers. Dordrecht-Boston-London.
\par
Walsh, J.B. (1986). An introduction to stochastic partial
differential equations, {\it Ecole d'Et\'e de Prob. de St.} Flour
XIV, 1984, Lect. Notes in Math 1180, Springer Verlag.
\par
Yong, J., and Zhou, X.Y. (1999). { Stochastic controls:
Hamiltonian systems and HJB equations}. New York: Springer-Verlag.
\par
Yosida, K. (1965).   {\it Functional Analysis}. Springer-Verlag.
Berlin, Gottingen, Heidelberg.
\par
 Zhou, X.Y. (1992). { A duality analysis on stochastic partial
differential equations}, {\it Journal of Functional Analysis} {\bf
103}, No. 2, 275--293.
\end{document}